\documentclass[12pt]{article}
\usepackage{amsmath,amssymb,amsthm,latexsym}

\newtheorem{cor}{Corollary}
\newtheorem{theorem}{Theorem}

\newcommand{\R}{\mathbb{R}}

\numberwithin{equation}{section}

\begin{document}

\title{ \textbf{A Note on the Ruin Problem\\
with Risky Investments}}
\author{David G Maher\\
School of Mathematics\\
University of New South Wales\\
Sydney  2052  Australia\\
Email: {\tt dmaher@maths.unsw.edu.au}}
\date{}
\maketitle

\begin{abstract}
{\small
We reprove a result concerning certain ruin in the classical problem of the probability of ruin with risky investments and several of its generalisations.  We also provide the combined transition density of the risk and investment processes in the diffusion case.}
\end{abstract}

{\small {\it Keywords:} Risk theory, ruin theory, Brownian motion, hyperbolic space.

}

\section{Introduction}

The classical Cram\'er-Lundberg risk model may be written as

\begin{equation}
X_t = u + ct - \sum_{k=1}^{N_t} Y_t
\end{equation}

where $X_t$ is the capital of the insurance company, $u$ is the initial capital, $c$ is the premium rate, $N_t$ is a Poisson process of rate $\lambda$, and $\{ Y_t \}_{k = 1}^{N_t}$ is a sequence of positive i.i.d. random variables, modelling the claim sizes (independently of $N_t$).\\

The probability of ruin is written as 
\begin{equation}
\Psi (u) = P \Bigl( X_t = u + ct - \sum_{k=1}^{N_t} Y_t \leq 0, \; \text{for some} \; t > 0 \Bigr)
\end{equation}

However, this model assumes no return on investments.  For many insurers, the extremely high competitiveness of today's financial market means that they actually have a zero or negative operating profit, thus relying on investing to make up the shortfall and make a profit.  The classical ruin model with risky investments has been considered by several authors, \cite{CAI}, \cite{CON}, \cite{FKP}, \cite{KN}, \cite{N}-\cite{P3} (and the references contained therein) being but a few.\\

In this paper, we shall firstly consider loss distributions with infinite support, reproving a result originating in \cite{P1}, and generalised in \cite{P2} and \cite{PZ}, that if the volatility of the investments is of a certain magnitude, then ruin is inevitable.  Our method of proof of this result is somewhat simpler and intuitively easier to understand.  This method is also quite flexible, and we are able to prove generalisations of this result for more general risk processes, as well as when the investment is modelled by a certain L\'evy process.\\

Additionally, we provide the transition density when the risk process is regarded as a diffusion, and conclude with some remarks on L\'evy processes.

\section{The main idea}

Consider the following: Take the classical Cram\'er-Lundberg risk model, and then act on the capital position by a geometric Brownian motion that models the investment:

$$
dX_t' = e^{\sigma B_t + \alpha t} \, dX_t
$$

That is,

\begin{equation}
X_t' = e^{\sigma B_t + \alpha t} u + \int_0^t e^{\sigma B_s + \alpha s} \, dX_s
\end{equation}

This is a version of risk process modified for investment considered in \cite{CAI}, \cite{KN} and \cite{P3}.  The continuous interaction of the combined (standard) risk process and investment can be described in the manner of a \textit{semi-direct product}:

$$
(X_t, e^{\sigma B_t + \alpha t})(X_s, e^{\sigma B_s + \alpha s}) = (X_t + X_s e^{\sigma B_t + \alpha t}, e^{\sigma B_{t+s} + \alpha (t+s)})
$$

The group with this operation is better known as \textit{real hyperbolic space}.  Some knowledge of hyperbolic space will be needed for Brownian motion, which we will be considering in section 3.  It is not a requisite for Theorem 1, although the idea of the action of investments on the risk process will quickly prove this result.

\subsection{Certain ruin with risky investments}

Consider the typical model for a share at time $t$, $S_t$,  expressed as geometric Brownian motion with drift and diffusion parameters of $a$ and $\sigma^2$, respectively.  That is,

\begin{equation}
dS_t = a S_t dt + \sigma S_t dB_t
\end{equation}

The solution of this S.D.E. is

$$ S_t = S_0 \exp \bigl( (a - \tfrac{\sigma^2}2) t + \sigma B_t \bigr) = S_0 \exp \bigl( \alpha t + \sigma B_t \bigr)$$

setting $\alpha = a - \tfrac{\sigma^2}2$.  The classical Cram\'er-Lundberg ruin model (1.1) is typically modified for investments by putting

\begin{equation}
X_t^\prime = u + a \int_0^t X_s^\prime ds + \sigma \int_0^t X_s^\prime dB_s + ct - \sum_{k=1}^{N_t} Y_t
\end{equation}

so that $X_t^\prime$ describes the evolution of the capital of an insurer which is continuously invested in an asset which follows a geometric Brownian motion (independent of $Y_t$ and $N_t$) with parameters $a$ and $\sigma$.  However, by using the model from (2.1):

\begin{equation}
X_t' = \int_0^t e^{\sigma B_s + \alpha s} \, dX_s
\end{equation}

we view the geometric Brownian motion (share price) as acting as a dilation on the classical risk model.  From this fact, our first result follows rather easily, proving the following:\\

\begin{theorem}
Consider the classical Cram\'er-Lundberg risk model with investments as in (2.4).  Assume further that the distribution of $Y_1$ does not have finite support, that is, $P(Y_1 > y) > 0$ for all $y > 0$.\\

If $\alpha < 0$ (or equivalently, $\frac{2a}{\sigma^2} < 1$), then ruin is certain.\\
\end{theorem}

\begin{proof}  It is enough to prove that the capital position $X_t^\prime$ is bounded, since the claim size can be large enough to ruin the company.  Now consider the classical risk model.  We have

\begin{equation}
X_t = u + ct - \sum_{k=1}^{N_t} Y_t \leq u + ct
\end{equation}

for all $t \geq 0$, since the $Y_t$ are positive distributions.  Thus from our model in (2.4),

\begin{align*}
X_t' & = e^{\sigma B_t + \alpha t}u + \int_0^t e^{\sigma B_s + \alpha s} \, dX_s\\
& \leq e^{\sigma B_t + \alpha t}u + \int_0^t e^{\sigma B_s+ \alpha s} \,d (cs)\\
& = e^{\sigma B_t + \alpha t}u + c \int_0^t e^{\sigma B_s+ \alpha s} \,ds
\end{align*}

We are interested in what happens for $t$ large.  We have

\begin{equation}
\exp \bigl( \sigma B_t + \alpha t \bigr) = \exp \bigl( t \, (\sigma B_t / t + \alpha) \bigr)
\end{equation}\\

By the strong law of large numbers, $B_t / t \rightarrow 0$ as $t \rightarrow \infty$, so our dilation term acts as $e^{t \alpha} \rightarrow 0$ as $t \rightarrow \infty$.  From this we deduce that $e^{\sigma B_t + \alpha t}u \rightarrow 0$ as $t \rightarrow 0$, and $\int_0^t e^{\sigma B_s+ \alpha s} \,ds$ is bounded for all $t > 0$.  Therefore, the capital position is bounded, and the theorem follows.

\end{proof}

Theorem 1 also holds for many generalisations of the Cram\'er-Lundberg ruin model.  The first concerns varying premium rates (c.f. \cite{PZ})

\begin{cor}  Suppose the premium rate in the risk model (1.1) is a bounded function, $c_t$.  Then with the assumptions of Theorem 1, ruin is certain.

\end{cor}

\begin{proof} This follows by putting $c = \sup_{t>0} c_t$ in (2.5) above.
\end{proof}

Another concerns when $N_t$ is a counting process other than the Poisson process (see, for example, \cite{DUF}).

\begin{cor}  Suppose $N_t$ is an arbitrary counting process in the risk model (1.1).  Then with the assumptions of Theorem 1, ruin is certain.

\end{cor}

We may also consider variations in the investment model where the Brownian motion is replaced by a more general L\'evy process:

\begin{cor}  Suppose the investments in the risk model (2.4) are modelled by the exponential functional

$$e^{\sigma L_t + \alpha t}$$

where $L_t$ is L\'evy process with mean 0.  Then with the assumptions of Theorem 1, ruin is certain.

\end{cor}

\begin{proof} This follows again by the strong law of large numbers since $\lim_{t \rightarrow \infty} L_t / t = 0$.
\end{proof}

{\bf Remark:}  These results also sheds some light on dividend constraints, since paying a dividend may be regarded as subtracting from the value of $a$ above, thus contributing to the overall probability of ruin.  Additionally, setting $\sigma = 0$ describes the risk model with a deterministic force of interest, $a$.

\section{The diffusion limit of the probability of ruin}

\subsection{The diffusion model}

The following characterisation of the probability of ruin was first introduced by Grandell in \cite{GRA}, who constructed a sequence of risk processes that converged weakly in the Skorohood topology to a diffusion process.\\

Put $C_t = \sum_{k=1}^{N_t} Y_t$, which has mean $\lambda \mu t$ and variance $\lambda m t$.  If the premium rate is set equal to $\lambda \mu$ then ruin will be certain.  To avoid this, a \textit{safety loading} $\rho$ is added as follows:

\begin{equation}
c = (1 + \rho) \lambda \mu \; \; \Longrightarrow \; \; \rho = \frac{c-\lambda \mu}{\lambda \mu}
\end{equation}

By regarding the risk process as a diffusion, we may use the classical result of the hitting times of Brownian motion to give the well known ``back of the envelope'' calculation of the probability of ruin by considering the \textit{diffusion limit of the probability of ruin}:

\begin{equation}
\Psi_D = \exp \Bigl\{ -2.\frac{\rho \mu u}{m} \Bigr\}
\end{equation}

We now considner the classical Cram\'er-Lundberg ruin model (considered as a diffusion), and then act on the capital position by a geometric Brownian motion that models the investment.  These models follow shifted and dilated Brownian motion on the groups $(\R, +)$ and $(\R^+, \times)$, respectively, with their interaction again described by the semi-direct product from (2.1).  That is, we will consider $\sum_{k=1}^{N_t} Y_t$ to be a Brownian motion with a drift of $\lambda \mu$ and diffusion co-efficient $\lambda m$.  This is somewhat at odds with \cite{GMY}, who argue that {\it no} financial asset may be correctly modelled on a continuous martingale.  However, some rigour for using the pure diffusion model has been provided in \cite{PG}, who give conditions for a weak convergence to a diffusion (ie, in the Skorohood topology as in \cite{GRA}), so that we may consider the ``unscaled version'' of the diffusion limit of the probability of ruin with investments.








\subsection{Real hyperbolic space}

Real hyperbolic space may be defined in several ways.  It is usually recognised as the \textit{Poincar\'e upper half-plane} (c.f. \cite{HE2}, and \cite{TER}, Ch III).  For our purposes, we identify real hyperbolic space as the semi-direct product $\R \rtimes \R^+$, where $\R$ is the group of reals under addition, and $\R^+$ is the group of positive reals under multiplication.  The semi-direct product $\R \rtimes \R^+$ is the topological space $\R \times \R^+$ with group multiplication given by

$$
(x, y)(x', y') = (x + x'y, yy')
$$

There are some slight technical details when performing analysis on this group.  Firstly, it is a prime example of a \textit{non-unimodular group}, ie, the left Haar measure is not equal to right Haar measure.  More importantly, it is {\it the only} simple, simply connected Lie group whose Laplacian cannot be written as a sum of squares of its vector fields (c.f. \cite{PAU}, Thm. 4.1).  Since the generator of Brownian motion is the Laplacian, this fact has implications for Brownian motion (c.f. \cite{PAU}, Cor. 4.4).  However, this is overcome in \cite{BOU} and \cite{COR} by considering a distinguished sub-Laplacian to generate a Brownian motion.\\

\subsection{Brownian motion on hyperbolic space}

Brownian motion on $\R \rtimes \R^+$ was studied explicitly in \cite{BOU} and \cite{COR}.  The random variable considered in \cite{BOU} was 

$$
\Bigl( \int_0^t e^{B_s} dW_s, \, e^{B_t} \Bigr)
$$

In \cite{COR} this was shown to be equivalent under the Dambis, Dubins-Schwarz theorem to the process

$$
(W_{A_t}, \, e^{B_t})
$$

where $(W_t)_{t \geq 0}$ and $(B_t)_{t \geq 0}$ are standard Brownian motions, and 

$$A_t = \int_0^t e^{2B_s}ds$$

Its characteristic function is 

\begin{align*}
E(\exp(i\xi W_{A_t} + i\zeta  e^{B_t})) & = E(\exp( i\zeta  e^{B_t}) (\exp(i\xi W_{A_t})| \mathcal{F}^B_t))\\
& = E(\exp( i\zeta  e^{B_t}) \exp(- \tfrac 12 \| \xi \|^2 A_t)) 
\end{align*}

This is then determined, and then inverted to obtain the density:

\begin{equation}
p_t(z, e^x) = \frac{e^{-x^2/2 t}}{\sqrt{2\pi t}} \int_0^\infty \frac{1}{\sqrt{2\pi y^2}} \exp \biggl(-\frac{z^2}{y^2} \biggr) a_t (x,y) dy
\end{equation}\\

However, we wish to consider $X_t$ as a Brownian motion on $\R$ starting from $u$, with mean $\rho \lambda \mu t$ and variance $\lambda m t$, being acted upon by a geometric Brownian motion on $\R^+$ with the parameters $\alpha$ and $\sigma$.  That is, we are considering

\begin{equation}
(X_{A_{\sigma^2 t}^{(\alpha)}}, \, e^{\sigma B_t + \alpha t})
\end{equation}

where $(X_t)_{t \geq 0}$ is as in (1.2), $(B_t)_{t \geq 0}$ is a standard Brownian motion, and 

\begin{equation}
A_{\sigma^2 t}^{(\alpha)} = \int_0^t e^{2 (\sigma B_s + \alpha s)}ds
\end{equation}

Equivalently from \cite{BOU} we have

$$\Bigl( \int_0^t e^{\sigma B_s + \alpha s} dX_s, \, e^{\sigma B_t + \alpha t} \Bigr)$$

which are the random variables  -  considered separately  -  in \cite{KN}.  Although as remarked in \cite{KN} that the density of $\int_0^t e^{\sigma B_s + \alpha s} dX_s$ cannot be put in a closed form, in the next section we shall derive the combined density of these processes.\\

The conditional distribution for $A_{\sigma^2 t}^{(\alpha)}$ in (3.5) comes from the functional studied in \cite{YOR}

$$A_t^{(\alpha)} = \int_0^t e^{2 (B_s + \alpha s)}ds$$

and write $A_t$ when $\alpha = 0$.  This functional originates in Yor's work with Bessel functions and its application to pricing Asian options.  For a given $B_t$, if the density of $A_t$ is given by

$$P(A_t \in du | \, B_t = x) = a_t(x,u) du$$

then

$$\frac{1}{\sqrt{2\pi t}} e^{-\frac{x^2}{2t}} a_t(x,u) = \frac{1}{u} \exp(-(1+e^{2x})/2u) \Theta_{\frac{1}{u}e^x} (t)$$

where

$$\Theta _{\frac{1}{u}e^x} (t) = \frac{x}{u \sqrt{2\pi^3 t}} \int_0^\infty e^{-\frac{y^2}{2t}} \exp(-e^x \cosh (y) / u) \sinh (y) \sin(\pi y/t) dy
$$

\newpage

\subsection{Transition density of the diffusion model}

Using the method described above, we are now in a position to obtain the transition density of the classical Cram\'er-Lundberg ruin model with investments, where both are considered to be diffusions.

\begin{theorem}
The transition density of the risk process (1.1) with investments

\begin{equation}
\Bigl( \int_0^t e^{\sigma B_s + \alpha s} dX_s, \, e^{\sigma B_t + \alpha t} \Bigr)
\end{equation}

is given by

\begin{equation}
p_t(z, e^x) = \frac{e^{-x^2/2\sigma^2 t}}{\sqrt{2\pi \sigma^2 t}} \int_0^\infty \frac{1}{\sqrt{2\pi \lambda \mu y^2}} \exp \biggl(-\frac{(z-\rho \lambda \mu t - u)^2}{\lambda \mu y^2} \biggr) a_{\sigma^2 t}(x,y) dy
\end{equation}

\end{theorem}

\begin{proof}  By the scalar invariance of $B_t$, and applying the Girsanov theorem, it is readily seen that the density of 

$$\int_0^t e^{2 (\sigma B_s + \alpha s)}ds = A_{\sigma^2 t}^{(\alpha)}$$

is given by

$$P(A_{\sigma^2 t}^{(\alpha)} \in du | \, \sigma B_t + \alpha t = x) = a_{\sigma^2 t}(x,u) du
$$\\

The characteristic function for the random variable in (3.6) is 

\begin{align*}
E \Bigl[\exp(i\xi X_{A_{\sigma^2 t}^{(\alpha)}} & + i\zeta  e^{(\sigma B_t + \alpha t)}) \Bigr] = E\Bigl[\exp( i\zeta  e^{(\sigma B_t + \alpha t)} ) (\exp(i\xi X_{A_{\sigma^2 t}^{(\alpha)}})| \mathcal{F}^B_t) \Bigr]\\
& = E\Bigl[\exp( i\zeta  e^{(\sigma B_t + \alpha t)} ) \exp( ( i \rho \lambda \mu t \xi + iu\xi - \tfrac 12 \lambda m \xi^2) A_{\sigma^2 t}^{(\alpha)}) \Bigr]\\
& = e^{i \rho \lambda \mu t \xi + iu \xi} \int_{-\infty}^\infty \exp(i\zeta e^{(\sigma x + \alpha t)}) \int_0^\infty \exp( -\tfrac 12 (\sqrt{\lambda m} \, \xi)^2)\\
& \phantom{abcdeabcdeabcdeabcdeabcdeabcde} \times a_{\sigma^2 t} (x,u) e^{-x^2/(2 \sigma^2 t)} du \, dx
\end{align*}

This is merely a shifted and dilated version of the characteristic function of the Brownian motion in \cite{COR}.  So, mutatis mutandis, we invert the Fourier transform to obtain the density:

$$ p_t(z, e^x) = \frac{e^{-x^2/2\sigma^2 t}}{\sqrt{2\pi \sigma^2 t}} \int_0^\infty \frac{1}{\sqrt{2\pi \lambda \mu y^2}} \exp \biggl(-\frac{(z-\rho \lambda \mu t - u)^2}{\lambda \mu y^2} \biggr) a_{\sigma^2 t}(x,y) dy $$\\

as required.

\end{proof}

Although the expression for $p_t(z, e^x)$ is far from tractable, we may use it to say the following:

\begin{cor}
Suppose $p_t(z, e^x)$ is the transition density for the risk process, then the probability that the company is ruined at time $t$ is given by
$$\int_{-\infty}^0 \int_{-\infty}^\infty p_t(z, e^x) dx \, dz$$
\end{cor}

Evaluating this integral will (more than likely) need to be done using numerical techniques.  This expression is particularly important in actuarial practice, since it may only be necessary for the insurer to be solvent at certain times, that is, review dates.\\

{\bf Remark:}  The above calculations would appear to be valid when the generator is that of the classical risk process.  That is, the process on $\R$ is a (discontinuous) semimartingale, rather than just a Brownian motion.  This allow us to compute the transition density in the case of a compound Poisson process, rather than just the diffusion model, where the combined density would be an expression similar to (3.7).  More generally, we could consider the risk process as a L\'evy process.  These models were discussed in detail in \cite{MS}, which we refer the reader to for many explicit examples.  The critical step is a ``Dambis, Dubins-Schwarz"-type theorem for discontinuous semimartingales, which is outside the scope of this paper.\\

{\bf Remark:}  Similarly, the transition density when the investment model is other than geometric Brownian motion requires significant modifications to Yor's work in generalising the density of the functional $A_t$ when the investment model is a L\'evy process rather than geometric Brownian motion.  That said, a functional of the form $A^\prime_t = \int_0^t e^{2L_s} ds$ where $L_t$ is a L\'evy process does provide us with a c\'adl\'ag modification to the risk model in (3.6).

\end{document}